\documentclass{article}

\usepackage[english]{babel}
\usepackage{makeidx}  

\usepackage{algorithm}
\usepackage{algpseudocode}
\usepackage{algorithmicx} 

\usepackage{amssymb}
\usepackage{amsmath}

\usepackage{graphicx}
\usepackage{booktabs}

\usepackage{import}
\usepackage{color}

\usepackage{caption}

\usepackage{titling}
\usepackage[numbers]{natbib}

\usepackage{geometry}
\newcommand{\X}{\mathcal{X}}
\newcommand{\Y}{\mathcal{Y}}
\newcommand{\A}{\mathbf{A}}

\DeclareMathOperator*{\argmin}{argmin}
\DeclareMathOperator*{\Ran}{Ran}

\begin{document}

\title{A warped kernel improving robustness in Bayesian optimization via random embeddings}

\author{Micka\"{e}l Binois$^{1,2}$ \and David Ginsbourger$^3$  \and Olivier Roustant$^1$}

\date{}

\maketitle             

\vspace{-1cm}

 \noindent $^1$ Mines Saint-\'{E}tienne, UMR CNRS 6158, LIMOS, F-42023 Saint-\'{E}tienne, France,\\ firstname.lastname@mines-stetienne.fr\\
 $^2$ Renault S.A.S., 78084 Guyancourt, France\\
$^3$ University of Bern, Department of Mathematics and Statistics, Alpeneggstrasse 22,
CH-3012 Bern, Switzerland, 
david.ginsbourger@stat.unibe.ch

\begin{abstract}
This works extends the Random Embedding Bayesian Optimization approach by integrating a warping of the high dimensional subspace within the covariance kernel.
The proposed warping, that relies on elementary geometric considerations, allows mitigating the drawbacks of the high extrinsic dimensionality while avoiding the algorithm to evaluate points giving redundant information. 
It also alleviates constraints on bound selection for the embedded domain, thus improving the robustness,
as illustrated with a test case with 25 variables and intrinsic dimension 6.\\

\end{abstract}

\section{Introduction}

The scope of Bayesian Optimization methods is usually limited to moderate-dimensional problems \cite{Koziel2011}. 
To overcome this restriction, \cite{Wang2013} recently proposed to extend the applicability of these methods to up to billions of variables, when only few of them are actually influential, through the so-called Random EMbedding Bayesian Optimization (REMBO) approach. In REMBO,  
optimization is conducted in a low-dimensional domain $\mathcal{Y}$, randomly embedded in the high-dimensional source space $\mathcal{X}$.
New points are chosen by maximizing the Expected Improvement (EI) criterion \cite{Mockus1978} with Gaussian process (GP) models incorporating the considered embeddings via two kinds of covariance kernels proposed in \cite{Wang2013}. 
A first one, $k_{\mathcal{X}}$, relies on Euclidean distances in $\mathcal{X}$. 
It delivers good performance in moderate dimension, albeit its main drawback is to remain high-dimensional so that the benefits of the method are limited. 
A second one, $k_{\mathcal{Y}}$, is defined directly over $\mathcal{Y}$ and is therefore independent from the dimension of $\X$. 
However, it has been shown \cite{Wang2013} to possess artifacts that may lead EI algorithms to spend many iterations exploring equivalent points.\\

Here we propose a new kernel with a warping (see e.g. \cite{Snoek2014}) inspired by simple geometrical ideas, that retains key advantages of $k_{\mathcal{X}}$ while remaining of low dimension like $k_{\Y}$.
Its effectiveness is illustrated on a 25-dimensional test problem with 6 effective variables.

\section{Background on the REMBO method and related issues}

The considered minimization problem is to find $\textbf{x}^* \in \argmin_{\textbf{x} \in \mathcal{X}} f(\textbf{x})$, with $f: \mathcal{X} \subseteq \mathbb{R}^D  \to \mathbb{R}$, where $\mathcal{X}$ is a compact subset of $\mathbb{R}^D$, assumed here to be $[-1,1]^D$ for simplicity. From \cite{Wang2013}, one main hypothesis about $f$ is that its effective dimensionality is $d_e < D$: there exists a linear subspace $\mathcal{T} \subset \mathbb{R}^D$ of dimension $d_e$ such that $f(\textbf{x}) = f(\textbf{x}_{\top} + \textbf{x}_{\perp}) = f(\textbf{x}_{\top})$, $\textbf{x}_{\top} \in \mathcal{T}$ and $\textbf{x}_{\perp} \in \mathcal{T}^{\perp} \subset \mathbb{R}^D$ (\cite{Wang2013}, Definition 1). Given a random matrix $\A \in \mathbb{R}^{D \times d}$ ($d \geq d_e$) with components sampled independently from $\mathcal{N}(0,1)$, for any optimizer $\textbf{x}^* \in \mathbb{R}^D$, there exists at least a point $\textbf{y}^* \in \mathbb{R}^d$ such that $f(\textbf{x}^*) = f(\textbf{Ay}^*)$ with probability 1 (\cite{Wang2013}, Theorem 2.). To respect box constraints, $f$ is evaluated at $p_{\mathcal{X}}(\textbf{A}\textbf{y})$, the convex projection of $\textbf{Ay}$ onto $\X$.
The low dimensional function to optimize is then $g: \mathbb{R}^d \to \mathbb{R}$, $g(\textbf{y}) = f\left(p_{\mathcal{X}}(\textbf{Ay})\right)$.\\

Optimizing $g$ is carried out using Bayesian Optimization, e.g, with the EGO algorithm \cite{Jones1998}. It bases on Gaussian Process Regression \cite{Rasmussen2006}, also known as Kriging \cite{Matheron1963}, to create a surrogate of $g$. Supposing that $g$ is a sample from a GP with known mean (zero here to simplify notations) and covariance kernel $k(.,.)$, conditioning it on $n$ observations $\mathbf{Z} = f(\textbf{x}_{1:n}) = g(\textbf{y}_{1:n}) $, provides a GP $Z(.)$ with mean $m(\textbf{x}) =  \textbf{k}(\textbf{x})^T K^{-1} \textbf{Z}$ and kernel $c(\textbf{x},\textbf{x}') = k(\textbf{x},\textbf{x}') - \textbf{k}(\textbf{x})^T K^{-1} \textbf{k}(\textbf{x}')$,    
where $\textbf{k}(\textbf{x}) = (k(\textbf{x}, \textbf{x}_i))_{1 \leq i \leq n}$ and $K = (k(\textbf{x}_i, \textbf{x}_j))_{1 \leq i,j \leq n}$. The choice of $k$ is preponderant, since it reflects a number of beliefs about the function at hand. Among the most commonly used are the ``squared exponential'' (SE) and ``Mat\'{e}rn" stationary kernels, with hyperparameters such as length scales or degree of smoothness \cite{Roustant2012,Stein1999}.
For REMBO, \cite{Wang2013} proposed two versions of the SE kernel with length scales $l$, namely the low-dimensional $k_{\mathcal{Y}}(\textbf{y}, \textbf{y}') = \exp \left(-\|\textbf{y} - \textbf{y}' \|^2_d \slash 2l_{\Y}^2\right)$ and the high-dimensional $k_{\mathcal{X}}(\textbf{y}, \textbf{y}')$ $=$ $\exp \left(-\left\|p_{\mathcal{X}}(\textbf{Ay}) - p_{\mathcal{X}}(\textbf{Ay}') \right\|^2_D \slash 2l_{\X}^2 \right)$ ($\textbf{y}, \textbf{y}' \in \Y$).\\ 

Selecting the domain $\mathcal{Y} \subset \mathbb{R}^d$ is a major difficulty of the method: if too small, the optimum may not be reachable while a too large domain renders optimizing harder, in particular since $p_{\X}$ is far from being injective. 
Distant points in $\mathcal{Y}$ may coincide in $\mathcal{X}$, especially far from the center, so that using $k_{\Y}$ leads to sample useless new points in $\Y$ corresponding to the same location in $\X$ after the convex projection.
On the other hand, $k_{\mathcal{X}}$ suffers from the curse of dimensionality when $\Y$ is large enough so that most or all of the points of $\X$ belonging to the convex projection of the subspace spanned by $\A$ onto $\X$ have at least one pre-image in $\Y$. 
Indeed, whereas embedded points $p_{\X}(\textbf{Ay})$ lie in a $d$ dimensional subspace when they are inside of $\X$, they belong to a $D$-dimensional domain when they are projected onto the faces and edges of $\X$.
To alleviate these shortcomings, after showing that with probability $1 - \epsilon$ the optimum is contained in the centered ball of radius $d_e / \epsilon$ (Theorem 3), the authors of \cite{Wang2013} then suggest to set $\Y = [-\sqrt{d}, \sqrt{d}]^d$. In practice, they split the evaluation budget over several random embeddings or set $d > d_e$ to increase the probability for the optimum to actually be inside $\mathcal{Y}$, slowing down the convergence. 

\section{Proposed kernel and experimental results}

Both $k_{\Y}$ and $k_{\X}$ suffering from limitations, it is desirable to have a kernel that retains as much as possible of the actual high dimensional distances between points while remaining of low dimension. This can be achieved by first projecting points orthogonally on the faces of the hypercube to the subspace spanned by $\A$: $\Ran(\A)$,
with $p_{\A}: \mathcal{X} \mapsto \mathbb{R}^D$, $p_{\A}(\textbf{x}) = \A(\A^T \A)^{-1} \A^T \textbf{x}$. Note that these back-projections from the hypercube can be outside of $\mathcal{X}$. 
The calculation of the projection matrix is done only once, inverting a $d \times d$ matrix. This solves the problem of adding already evaluated points: their back-projections coincide. 
Nevertheless, distant points on the sides of $\X$ from the convex projection can be back-projected close to each other, which may cause troubles with the stationary kernels classically used.\\

\begin{figure}[htpb]%
	\centering
		\def\svgwidth{0.9\textwidth}
	\begingroup%
  \makeatletter%
  \providecommand\color[2][]{%
    \errmessage{(Inkscape) Color is used for the text in Inkscape, but the package 'color.sty' is not loaded}%
    \renewcommand\color[2][]{}%
  }%
  \providecommand\transparent[1]{%
    \errmessage{(Inkscape) Transparency is used (non-zero) for the text in Inkscape, but the package 'transparent.sty' is not loaded}%
    \renewcommand\transparent[1]{}%
  }%
  \providecommand\rotatebox[2]{#2}%
  \ifx\svgwidth\undefined%
    \setlength{\unitlength}{616.63603516bp}%
    \ifx\svgscale\undefined%
      \relax%
    \else%
      \setlength{\unitlength}{\unitlength * \real{\svgscale}}%
    \fi%
  \else%
    \setlength{\unitlength}{\svgwidth}%
  \fi%
  \global\let\svgwidth\undefined%
  \global\let\svgscale\undefined%
  \makeatother%
  \begin{picture}(1,0.46749837)%
    \put(0,0){\includegraphics[width=\unitlength]{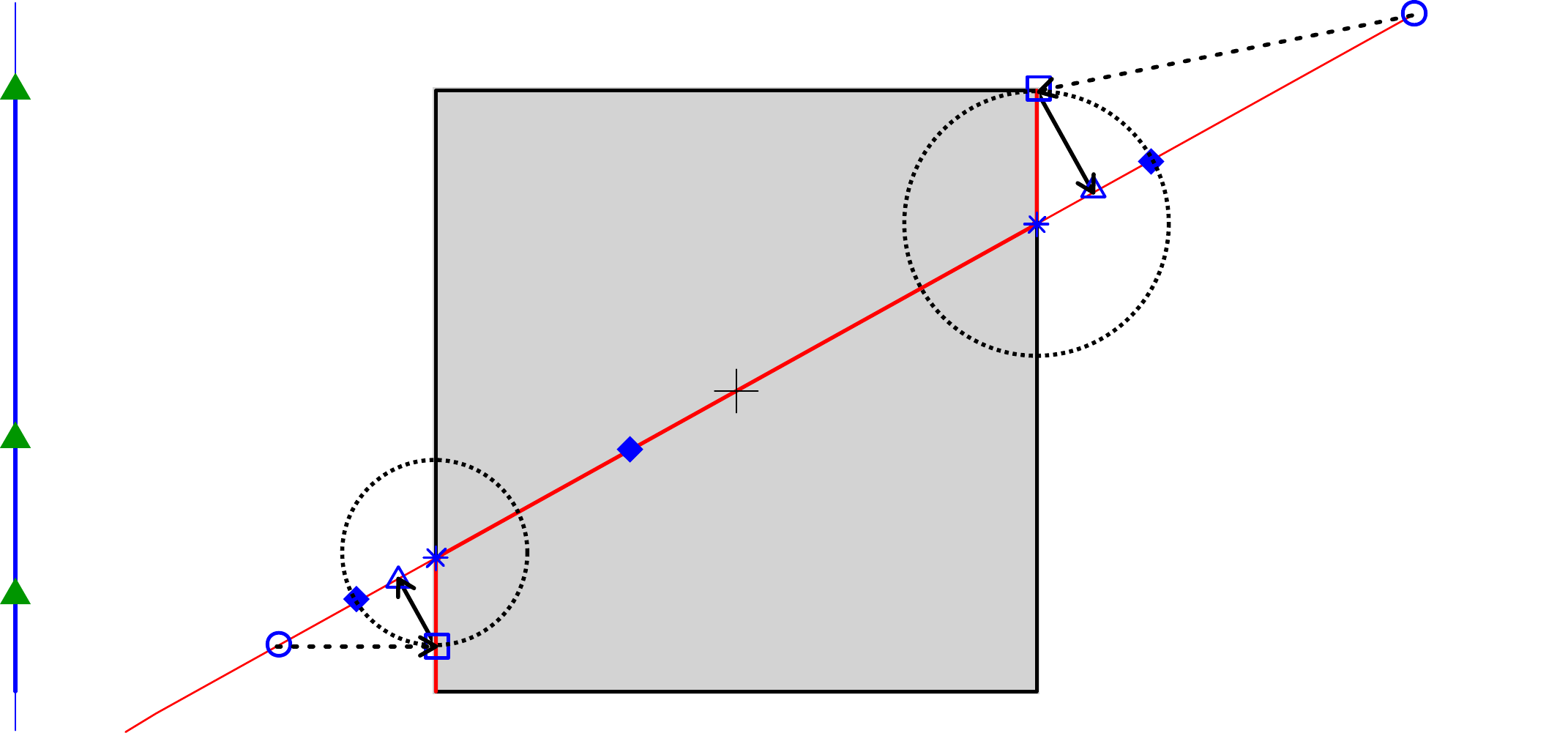}}%
    \put(0.02240199,0.40644494){\color[rgb]{0,0,0}\makebox(0,0)[lb]{\smash{$y_1$}}}%
    \put(0.02760739,0.18496767){\color[rgb]{0,0,0}\makebox(0,0)[lb]{\smash{$y_2$}}}%
    \put(0.02501266,0.08117834){\color[rgb]{0,0,0}\makebox(0,0)[lb]{\smash{$y_3$}}}%
    \put(0.45215871,0.24069216){\color[rgb]{0,0,0}\makebox(0,0)[lb]{\smash{$O$}}}%
    \put(0.37990002,0.13692658){\color[rgb]{0,0,0}\makebox(0,0)[lb]{\smash{$\Psi(y_2) = \textbf{A}y_2$}}}%
    \put(0.76036467,0.35822375){\color[rgb]{0,0,0}\makebox(0,0)[lb]{\smash{$\textbf{z}_1'' = \Psi(y_1)$}}}%
    \put(0.13796391,0.09808266){\color[rgb]{0,0,0}\makebox(0,0)[lb]{\smash{$\Psi(y_3)$}}}%
    \put(0.14198849,0.02227296){\color[rgb]{0,0,0}\makebox(0,0)[lb]{\smash{$\textbf{A}y_3$}}}%
    \put(0.87502636,0.41604851){\color[rgb]{0,0,0}\makebox(0,0)[lb]{\smash{$\textbf{A}y_1$}}}%
    \put(0.01964646,0.44815225){\color[rgb]{0,0,0}\makebox(0,0)[lb]{\smash{$\mathcal{Y}$}}}%
    \put(0.30421639,0.37037696){\color[rgb]{0,0,0}\makebox(0,0)[lb]{\smash{$\mathcal{X}$}}}%
    \put(0.69055457,0.31827716){\color[rgb]{0,0,0}\makebox(0,0)[lb]{\smash{$\textbf{z}_1$}}}%
    \put(0.60744302,0.3332298){\color[rgb]{0,0,0}\makebox(0,0)[lb]{\smash{$\textbf{z}'_1$}}}%
    \put(0.54862128,0.43078641){\color[rgb]{0,0,0}\makebox(0,0)[lb]{\smash{$p_{\mathcal{X}}(\textbf{A}y_1)$}}}%
  \end{picture}%
\endgroup%
	\caption{
	Illustration of the new warping $\Psi$ , $d = 1$ and $D = 2$, from triangles in $\Y$ to diamonds in $\X$, on three points $y_1, y_2,
y_3$. As for REMBO, the points $y_i$ are first mapped by $\A$ and convexly projected onto $\X$ (if out of $\X$). If the resulting image is
strictly contained in $\X$ -- as for $y_2$ -- nothing else is done.
Otherwise, the new warping is defined in two supplementary steps: back-projection onto Ran($\A$) (giving $\textbf{z}_i$) and stretching out in the
resulting line $[0, \textbf{z}_i)$ (red solid line) by reporting the distance between the intersection of $[0, \textbf{z}_i]$ on the frontier of $\X$, $\textbf{z}'_i$, and the initial convex projection $p_\X(\A y_i)$. The points $y_1$ and $y_3$
correspond to cases where such projections are on a corner or a face of $\X$.
	}%
	\label{fig:schema}%
\end{figure}
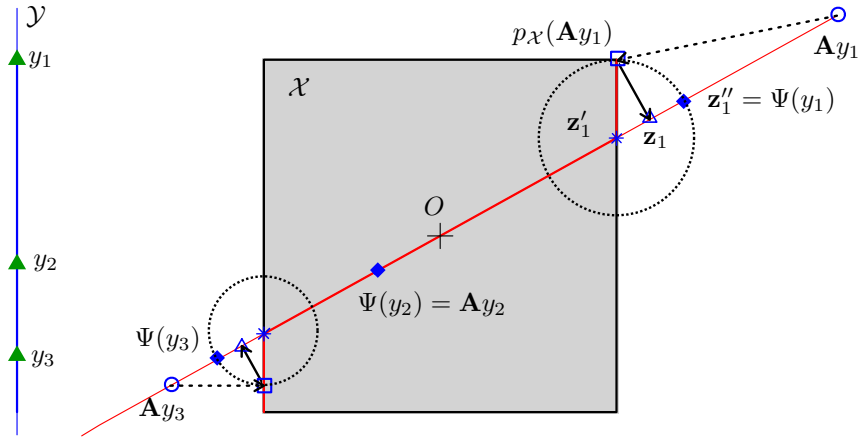

The next step is to respect as much as possible distances on the border of $\X$, denoted $\partial{\X}$. Unfolding and parametrizing the manifold corresponding to 
the convex projection of the embedding of $\mathcal{Y}$ with $\A$ would be best but unfortunately it seems intractable with high $D$. Indeed, it amounts to finding each intersection of the $d$-dimensional subspace spanned by $\A$ with the faces of the $D$-hypercube, before describing the parts resulting from the convex projection. 
Alternatively, we propose to distort the back-projections which are outside of $\X$,
corresponding to those convex-projected parts on the sides of $\partial{\X}$.
In more details, from the back-projection of the initial mapping with $p_{\X}$, a pivot point is selected as the  intersection between $\partial{\X}$ and the line ($O;p_{\A}(p_{\mathcal{X}}(\textbf{A}y))$. Then the back-projection is stretched out such that the distance between the pivot point and the initial convex projection are equal. It results in respecting the distance \emph{on the embedding} between the center $O$ and the initial convex projection.  
The resulting warping, denoted $\Psi$, is detailed in Algorithm \ref{alg:Psi} and illustrated in Figure \ref{fig:schema}. Based on this, any positive definite kernel $k$ on $\mathcal{Y}$ can be used. For example, the resulting SE kernel is $k_{\Psi}(\textbf{y}, \textbf{y'}) = \exp \left(-\left\|\Psi(\textbf{y}) - \Psi(\textbf{y}') \right\|^2_{D} \slash 2l_{\Psi}^2\right)$. Note that the function value corresponding to $\Psi(\textbf{y})$ remains $g(\textbf{y})$.\\

\begin{algorithm}[htpb]
\caption{Calculation of $\Psi$.}
\label{alg:Psi}
\begin{algorithmic}[1]
\State  Map $\textbf{y} \in \mathcal{Y}$ to $\textbf{Ay}$
\State \textbf{If} $\textbf{Ay} \in \mathcal{X}$ \textbf{Then} 
\State \hspace{0.5cm}Define $\Psi(\textbf{y) = \textbf{Ay}}$
\State \textbf{Else} 
\State \hspace{0.5cm} Project onto $\mathcal{X}$ and back-project onto $\Ran(\A)$: $\textbf{z} = p_{\textbf{A}}(p_{\mathcal{X}}(\textbf{Ay}))$
\State  \hspace{0.5cm} Compute the intersection of $[O; \textbf{z}]$ with $\partial{\X}$: $\textbf{z}' = (\max_{i=1, \dots, D} |z_i|)^{-1} \textbf{z}$
\State \hspace{0.5cm} Define $\Psi(\textbf{y}) = \textbf{z}' + \|p_{\mathcal{X}}(\textbf{A}y) - \textbf{z}' \|_D . \frac{\textbf{z}'}{\|\textbf{z}' \|_D}$
\State \textbf{EndIf}
\end{algorithmic}
\end{algorithm}

Like $k_{\X}$, $k_{\Psi}$ is not hindered by the non-injectivity brought by the convex projection $p_{\X}$.
Furthermore, it can explore sides of the hypercube without spending too much budget since belonging to $\Ran(\textbf{A})$ (all distances between embedded points after warping are $d$-dimensional instead of $D$-dimensional, thus smaller, hence limiting the risk of over-exploring sides of $\X$). 
It is thus possible to extend the size of $\mathcal{Y}$ to avoid the risk of missing the optimum. 
For instance, one can check that $\Y$ is larger than $[-\gamma, \gamma]^d$ with $\gamma$ such that $\gamma^{-1} = \min \limits_{j \in [1, \dots, D]} \sum\limits_{i = 1}^d {|A_{j,i}|}$, with $A_{j,i}$ the components of $\A$, ensuring to span $[-1,1]$ for each of the $D$ variables.\\

We compare the performances of the usual REMBO method with $k_{\mathcal{Y}}$, $k_\X$ and the proposed $k_{\Psi}$, with a unique embedding. Tests are conducted with the \emph{DiceKriging} and \emph{DiceOptim} packages \cite{Roustant2012}. We use the isotropic Mat\'{e}rn 5/2 kernel with hyperparameters estimated with Maximum Likelihood and we start optimization with space filling designs of size $10d$. Initial designs are modified such that no points are repeated in $\mathcal{X}$ for $k_{\mathcal{Y}}$ and $k_\X$. For $k_{\Psi}$, we apply $\Psi$ to bigger initial designs before selecting the right number of points, as distant as possible between each other. Experiments are repeated fifty times, taking the same random embeddings for all kernels. To allow a fair comparison, $\Y$ is set to $[-\sqrt{d},\sqrt{d}]^d$ for all kernels and the computational efforts on the maximization of the Expected Improvement are the same.\\
 
\begin{figure}[htpb]%
	\centering
	\def\svgwidth{0.8\textwidth}
	\begingroup%
  \makeatletter%
  \providecommand\color[2][]{%
    \errmessage{(Inkscape) Color is used for the text in Inkscape, but the package 'color.sty' is not loaded}%
    \renewcommand\color[2][]{}%
  }%
  \providecommand\transparent[1]{%
    \errmessage{(Inkscape) Transparency is used (non-zero) for the text in Inkscape, but the package 'transparent.sty' is not loaded}%
    \renewcommand\transparent[1]{}%
  }%
  \providecommand\rotatebox[2]{#2}%
  \ifx\svgwidth\undefined%
    \setlength{\unitlength}{368.53295898bp}%
    \ifx\svgscale\undefined%
      \relax%
    \else%
      \setlength{\unitlength}{\unitlength * \real{\svgscale}}%
    \fi%
  \else%
    \setlength{\unitlength}{\svgwidth}%
  \fi%
  \global\let\svgwidth\undefined%
  \global\let\svgscale\undefined%
  \makeatother%
  \begin{picture}(1,0.47555977)%
    \put(0,0){\includegraphics[width=\unitlength]{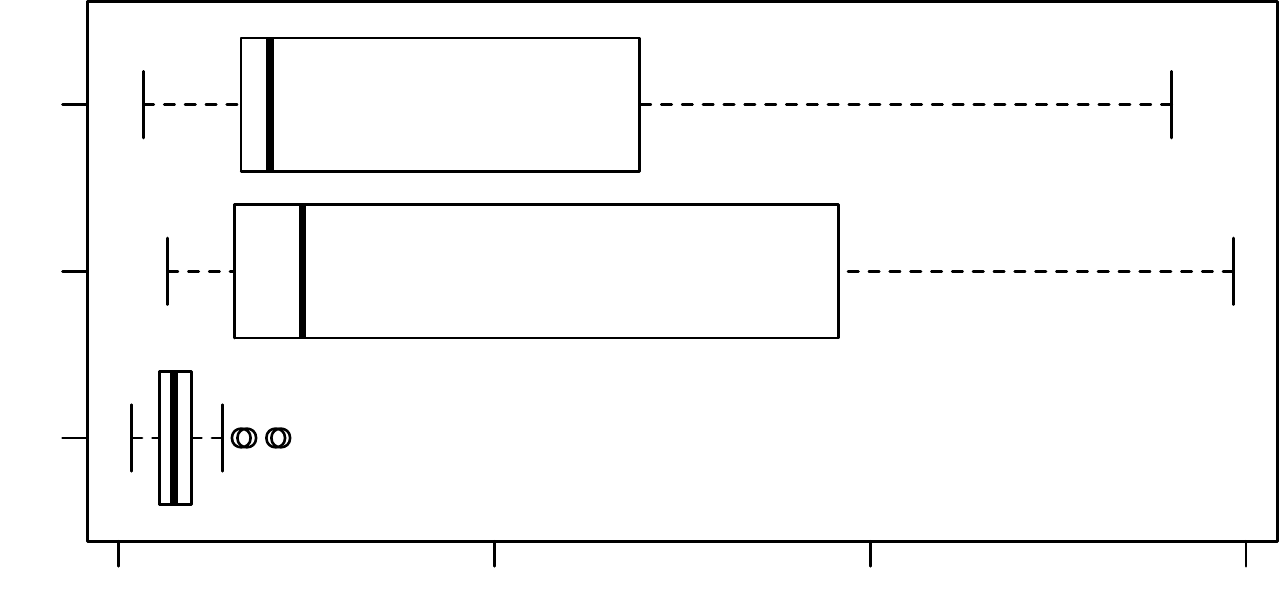}}%
    \put(0.06621467,0.00067207){\color[rgb]{0,0,0}\makebox(0,0)[lb]{\smash{0.0}}}%
    \put(0.36396186,0.00161542){\color[rgb]{0,0,0}\makebox(0,0)[lb]{\smash{0.5}}}%
    \put(0.65788263,0.00053005){\color[rgb]{0,0,0}\makebox(0,0)[lb]{\smash{1.0}}}%
    \put(0.95006814,0.00052997){\color[rgb]{0,0,0}\makebox(0,0)[lb]{\smash{1.5}}}%
    \put(-0.00185702,0.38386339){\color[rgb]{0,0,0}\makebox(0,0)[lb]{\smash{$k_{\mathcal{X}}$}}}%
    \put(-0.00185702,0.25393046){\color[rgb]{0,0,0}\makebox(0,0)[lb]{\smash{$k_{\mathcal{Y}}$}}}%
    \put(-0.00185702,0.12397126){\color[rgb]{0,0,0}\makebox(0,0)[lb]{\smash{$k_{\Psi}$}}}%
  \end{picture}%
\endgroup%
	\caption{Boxplot of the optimality gap (best value found minus actual minimum) for kernels $k_\X$, $k_{\Y}$ and $k_{\Psi}$ on the Hartmann6 test function (see e.g. \cite{Jones1998}) with $250$ evaluations, $d = d_e = 6$, $D= 25$.}%
	\label{fig:2}%
\end{figure}
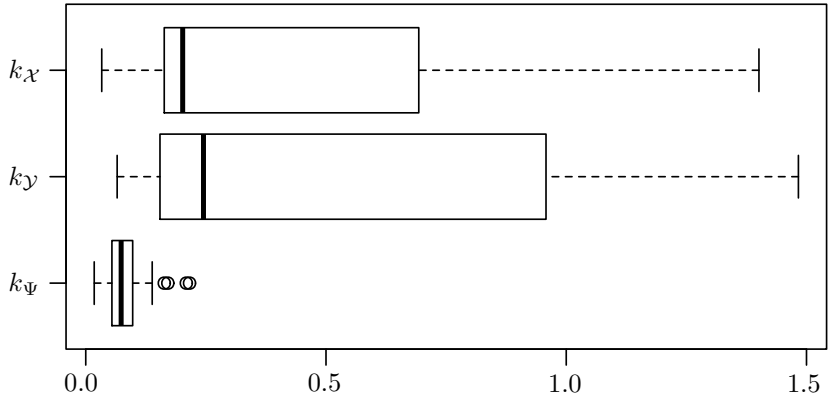

Results in Figure \ref{fig:2} show that the proposed kernel $k_{\Psi}$ outperforms both $k_{\Y}$ and $k_\X$ when $d = 6$. 
In particular, $k_\Y$ loses many evaluations on the sides of $\Y$ for already known points in $\X$ and $k_\X$ has a propensity to explore sides of $\X$, while $k_\Psi$ avoids both pitfalls.\\

\section{Conclusion and perspectives}

The composition with a warping of the covariance kernel used with REMBO wipes out some of the previous shortcomings. It thus achieved the goal of improving the results with a single embedding, as was shown on the Hartman6 example.
Studying the efficiency of splitting the evaluation budget between several random embeddings, compared to relying on a single one along with $k_{\Psi}$, would be the scope of future research. Of interest is also the study of the embedding itself, such as properties ensuring fast convergence in practice.\\

\section*{Acknowledgments}
This work has been conducted within the frame of the
ReDice Consortium, gathering industrial (CEA, EDF, IFPEN, IRSN, Renault) and academic (Ecole des Mines de Saint-Etienne, INRIA, and the University of Bern) partners around advanced methods for Computer Experiments.\\
The authors also thanks the anonymous reviewers as well as Frank Hutter for their helpful suggestions.

{\normalsize
\bibliographystyle{apalike}
\bibliography{biblio}}

\end{document}